\documentclass[12pt]{iopart}
\usepackage{url}
\usepackage{fancyvrb}
\usepackage{graphicx}
\usepackage[square,numbers]{natbib}
\usepackage{setstack}
\usepackage{iopams}

\eqnobysec

\begin{document}

%
%
%

\title[Construction of ODE systems]
      {Construction of ODE systems from time series data 
       by a highly flexible modelling approach}

\author{T Dierkes} 

\address{Konrad-Zuse-Zentrum f{\"u}r Informationstechnik Berlin (ZIB)\\ 
         Takustra{\ss}e 7, 14195 Berlin-Dahlem, Germany}

\ead{dierkes@zib.de}

\begin{abstract}
In this paper, a down-to-earth approach to purely data-based modelling of unknown 
dynamical systems is presented.  Starting from a classical, explicit ODE formulation 
$y' = f(t,y)$ of a dynamical system, a method determining the unknown right-hand 
side $f(t,y)$ from some trajectory data $y_{k}(t_{j})$, possibly very sparse, is given.  
As illustrative examples, a semi-standard predator-prey model is reconstructed from 
a data set describing the population numbers of hares and lynxes over a period of 
twenty years \citep{Dd2015_GUIDE}, and a simple damped pendulum system with a highly 
non-linear right-hand side is recovered from some artificial but very sparse data 
\citep{perona2000trajectory}.
%
\end{abstract}



%
%
\section{Introduction}

Modelling of dynamical systems is one of the key topics not only in systems biology, 
but also in other disciplines \citep{walter1997identification,Dd_NUM2}.  Most 
often, a desired behaviour of a dynamical system or, at least, parts of it, 
has already been anticipated.  In this paper we will present a down-to-earth 
approach to determine a completely unknown dynamical system from given, possibly 
very sparse, data sets; yet open enough for easy incorporation of additional
knowledge about the dynamical system to be considered.  Previous work related 
to ours include \citep{perona2000trajectory,eisenhammer1991modeling} and 
\citep{raue2015data}, describing a software package called \texttt{Data2Dynamics}.  
In the former two papers, a trajectory method is presented that first replaces 
the right-hand side of the ODE by a linear combination of a family of known 
functions, and tries then to find the coefficients of the linear combination 
by fitting the trajectory solutions to the given data.  The title of the third 
paper, ``Data2Dynamics : a modeling environment ...'', is 
a classical misnomer:  The complete dynamical model has already to 
be known in terms of an ODE system, only a possible control function for the ODE 
appears in their method as variable and unknown and could, in principle, be 
estimated, simultaneously with the classical fit approach, from the given data.  
In contrast to all three cited papers, our approach is slightly different, yet 
quite robust and fast and, therefore, extremely versatile, as will be explained 
in this paper.

\medskip

The dynamical system is supposed to be given as an explicit ODE,
\[
   y' = f(t,y) 
\]
where the dash $'$ represents the derivative w.r.t.~the time variable $t \in \mathbb{R}$ 
of the solution trajectory, $y = y(t)$, $y \in \mathbb{R}^{n}$.  The right-hand side, 
$f : \mathbb{R} \times \mathbb{R}^{n} \longrightarrow \mathbb{R}^{n}$, presumably 
smooth enough, has to be determined from data sampled of trajectory solutions, 
$y_{k}(t_{k,j})$, to this ODE, for $k \in \{1,\ldots,n\}$  and for some discrete 
time points $t_{k,1} < t_{k,2} < t_{k,3} < \ldots < t_{k,m}$.  Note that the sample 
time points can totally be asynchronous.  

\medskip

Since any non-autonomous system is readily transformed to an equivalent autonomous 
system, w.l.o.g.~we can assume $f$ to be independent of the time variable $t$, 
\[
   f(t,y) = f(y) \, . 
\]

If a trajectory component, $y_{k}$, has no sample points, we have to stipulate 
some values for that component describing a desired behaviour of the dynamical system 
in this component.  The general idea we have in mind here, is that only a subsystem 
of a \emph{known}, but apart from that, parametrised dynamical system, 
\[ 
   y' = f(y \, ; \, p) 
\]
with a fixed parameter vector $p \in \mathbb{R}^{q}$ is unknown, i.e.~only some of 
the solution trajectory components are given by sampled measurement points, or else 
by values describing a new, desired evolution of the dynamical system under investigation.

\section{Modelling Approach}

Having $m \in \mathbb{N}$ sample points of each solution trajectory component 
$y_{k,j} := y_{k}(t_{k,j})$, i.e.~$n m$ data values, of the unknown ODE system
\begin{equation} \label{eq:unknown_ode_system}
   y' = f(y)
\end{equation}
the basic idea of our approach to find $f$ from the given data is to first construct 
a suitable approximation $y^*(t)$ to the data points $y_{k,j}$ for that the derivative 
can readily be computed as well.  Then, we can construct a (discrete) mapping
\[
    F : y^* \mapsto (y^*)'
\]
with arbitrary many samples over the interval $[t_{min}, t_{max}]$, determined by the 
range of the sample time points of the given data.  

\medskip

From this, by making use of a set of ansatz functions 
$\{ \varphi_{\ell}\,,\ \ell = 0,1,\ldots \}$ over $\mathbb{R}^n$ (or $[a,b]^n$, 
if bounds of the components $y_k$ are available), it is possible to reconstruct 
the unknown function $f : \mathbb{R}^n \longrightarrow \mathbb{R}^n$ by linear,
or even non-linear, combinations of the $\varphi_{\ell}$: For as many 
$t_j \in [t_{min},t_{max}]$ as necessary, in a linear case, we eventually obtain
\begin{equation} \label{eq:lin_combi_ansatz}
    (y_{k}^{*})'(t_j) = \sum_{\ell = 0}^{L} \, 
                         c^{k}_{\ell} \,\, \varphi_{\ell}\left(\,y^{*}(t_j)\,\right), 
                         \quad
    k = 1,\ldots,n
\end{equation}
where, for each $k \in \{1,\ldots,n\}$, the $c^{k}_{\ell},\ \ell = 0, \ldots, L$, are 
coefficients to be determined in the least-squares sense.  Note that the right-hand
side depends on all $n$ trajectory components whereas the left-hand side includes only 
the $k$-th component.

\medskip

The solution of (\ref{eq:lin_combi_ansatz}) finally yields the description of the 
unknown $f$ in terms of the (deliberately chosen) ansatz functions $\varphi_{\ell}(\cdot)$.

\subsection{Best Polynomial Approximation}

Considering the set
\[
    \mathcal{P}_{N} = \{ \, \mathrm{polynomials\ on\ } [-1,1] 
                                   \mathrm{\ of\ degree\ at\ most\ } N \, \}
\]
it is well-known that for a function $g(t)$, Lipschitz continuous on $[-1,1]$, 
it is very difficult to find its corresponding \emph{minimax polynomial}, i.e.~the 
polynomial $p^*[g](t)$ that has the smallest maximum deviation from $g(t)$, 
see \citep{press2002_NR},
\[
   \|g - p^*[g]\|_{\infty} = \min\limits_{p \in \mathcal{P}_{N}} \|g - p\|_{\infty} \, .
\]
As it turns out, most remarkably, the corresponding Chebyshev series to $g$, 
which is absolutely and uniformly convergent,
\[
   g(t) =\sum_{k = 0}^{\infty} \, a_{k} \, T_{k}(t)   
\]
yields very nearly the same polynomial of degree $N$ as the \emph{minimax polynomial}
$p^*[g]$ when truncated at the $N$-th term,
\[
   p^*[g](t) \approx \sum_{k=0}^{N} \, a_{k} \, T_{k}(t) \, .
\]
Here, $T_{k}(t)$ denotes the $k$-th Chebyshev polynomial, defined as real part of
the complex function $z^{k}$ on the unit circle (see, e.g., \citep{trefethen2013_AT}),
\[
   T_{k}(t) = \mathrm{Re\ } z^{k} = 
              \frac{1}{2}(z^{k} + \overline{z}^{\,k}) = \cos(k\theta), \quad
   t = \cos \theta, \quad \mathrm{and\ \ \ } |z| = 1,
\]
and the coefficients $a_{k}$ are given for $k \geq 1$ by the formula
\begin{equation} \label{eq:cheb_coeffs}
   a_{k} = \frac{2}{\pi} \int\limits_{-1}^{1} \frac{g(t)\,T_{k}(t)}{\sqrt{1 - t^2}} \, dt
\end{equation}
and for $k = 0$ by the very same formula with $2/\pi$ changed to $1/\pi$.

\medskip

The discrete analogue to (\ref{eq:cheb_coeffs}), when $g(t)$ is only given at the 
$M+1$ Chebyshev nodes 
\[
   t_j = \cos\,\theta_j\, , \quad  \theta_j = \frac{\pi(j+1/2)}{M+1} \, , \quad j=0,\ldots,M \, , 
\]
reads \citep{quarteroni2007_NM, press2002_NR}
\begin{eqnarray} \label{eq:discrete_cheb_coeffs}
   \tilde{a}_{k} & = \frac{2}{M+1} \sum_{j=0}^{M} \,\, g(t_{j})\,T_{k}(t_{j}) \\
                 & = \frac{2}{M+1} \sum_{j=0}^{M} \,\, 
                                     g\left( \cos \frac{\pi(j+1/2)}{M+1} \right) \,
                                                   \cos \frac{\pi k(j+1/2)}{M+1} \, .
\end{eqnarray}

\medskip

In the very likely case that the sample time points $t_j$ of the given data do \emph{not} 
coincide with the Chebyshev nodes $\cos\,\theta_j$, we can use a simple linear 
interpolation scheme to transform the given data values to the required Chebyshev nodes.
An even more simple table lookup procedure would in some cases do equally well.
In the case of very rough data, alternatively, we could think of an additional smoothing
step by making use of a suitable spline interpolation where the data points are the 
control points of the spline interpolation scheme, for example \citep{de_boor1978_Guide}.  
Having such a smooth interpolation of the given data at hand, values at the desired 
Chebyshev nodes are then readily available.

\medskip

Thus, as building block for the first step, we arrive at the discrete Chebyshev approximation 
to $g$ of order $M$,
\begin{equation} \label{eq:discrete_approx}
   g(t) \approx \frac{1}{2} \tilde{a}_0 + \sum_{k=1}^{M} \, \tilde{a}_{k} \, T_{k}(t) \, ,
\end{equation}
evidently a polynomial of degree $M$.  Moreover, the derivative of $g(t)$ can readily 
be approximated in turn,
\begin{equation} \label{eq:discrete_deriv}
   g'(t) \approx \frac{1}{2} \tilde{a}_0' + \sum_{k=1}^{M} \, \tilde{a}_{k}' \, T_{k}(t) \, ,
\end{equation}
where the new coefficients $\tilde{a}_{k}'$ are given by the simple backwards recurrence 
formula
\begin{eqnarray*}
   \tilde{a}_{M+1}' & = \tilde{a}_{M}' = 0 \, , \\
   \tilde{a}_{k-1}' & = \tilde{a}_{k+1}' + 2\,k\,\tilde{a}_{k} \, ,\quad k = M, M-1,\ldots,1\,. 
\end{eqnarray*}

Now, making use of (\ref{eq:discrete_approx}) and (\ref{eq:discrete_deriv}), we 
can readily set up the discrete mapping $F : y^* \mapsto (y^*)'$ from the data $y_{k,j}$ 
by generating one separate Chebyshev approximation for each component $y_k$.  Depending
on the smoothness of the given data, the maximal degree $M$ of these approximations can
be adjusted accordingly that, as an intriguing side effect, comes in handy as a cheap 
filtering device in this step.

\subsection{Least-squares Solution: Linear and Non-linear Case}

For each $k \in \{1,\ldots,n\}$, equation (\ref{eq:lin_combi_ansatz}) 
can be rewritten as
\begin{equation} \label{eq:lsq_linear_system}
   b^{(k)} = A \, c^{(k)}
\end{equation}
where the $m$-vector $b^{(k)}$ and the $\bigg(m \times (L+1)\bigg)$-matrix $A$ 
(sometimes known as \emph{Gram matrix}) are given respectively by
\begin{eqnarray*}
   b_{j}^{(k)} & := \big(y_{k}^*\big)'(t_j) \, , & j = 1,\ldots,m \\
    A_{j,\ell} & := \varphi_{\ell}\bigg( \, y_{1}^*(t_j),\ldots,y_{n}^*(t_j) \, \bigg) \, , 
                   \quad & \ell = 0,\ldots,L  \, .
\end{eqnarray*}
Here, as explained in the previous section, the approximations $y_{k}^*(\cdot)$ and 
its derivatives are evaluated in terms of their truncated Chebyshev series of order $M$,
\[
   y_{k}^*(t) = \left[ \,\, \sum_{j=0}^{M} \tilde{y}_{k,j} \, T_j(t) \, \right] 
                                           - \frac{1}{2} \tilde{y}_{k,0}
\]
with the discrete Chebyshev coefficients $\tilde{y}_{k,j}$ computed by 
(\ref{eq:discrete_cheb_coeffs}) with the given $(n \times m)$-data matrix $(y_{k,j})$ .

\medskip

The solution $\left(c^{(k)}\right)^*$, in a least-squares sense, to the linear 
system (\ref{eq:lsq_linear_system}) is readily obtained by making use of the 
QR-decomposition $A = Q R P^{T}$ with pivoting, $P$ being a permutation matrix,
\begin{eqnarray*}
     \left(c^{(k)}\right)^* = P \, \tilde{R}^{-1} \, \tilde{Q}^{T} \, b^{(k)} \, ,
        \quad & R = \left[ \begin{array}{c} \tilde{R} \\ 0 \end{array} \right] \, , 
        \quad & Q = \left[ \begin{array}{cc} \tilde{Q} & \divideontimes \end{array} \right] \, , \\
     & \tilde{R} \in \mathbb{R}^{(L+1) \times (L+1)} \, , \quad 
     & \tilde{Q} \in \mathbb{R}^{m \times (L+1)} \, ,
\end{eqnarray*}
provided the \emph{Gram matrix} $A$ has full rank, and $m \geq L+1$.  
If one or both of these conditions are violated, there are adequate 
formulae available as well, see \citep{Dd_NUM1}.

\medskip

In the not so unlikely case, some (or all) of the unknown coefficients $c^{k}$ 
enter \emph{non-linearly} the ansatz function family $\{\, \varphi_\ell \, \}$, 
we have to resort to an iterative scheme, such as Gauss-Newton, in order to get 
a decent estimate of $\left(c^{k}\right)^*$.  
For the details we refer to \citep{Dd_newton}.

\subsection{Verification}

Up to now, we have not touched any of the dynamical properties of the underlying,
but unknown ODE system (\ref{eq:unknown_ode_system}), with the only exception of
(\ref{eq:discrete_deriv}) in determining the left-hand side $b^{(k)}$ in 
(\ref{eq:lsq_linear_system}).  Consequently, we can take advantage of this 
fact in the following manner if we want to verify the resulting solution of 
the previous section.

\medskip

With the least-squares solutions $(c^{(k)})^*$, $k=1,\ldots,n$, at hand, it 
is possible to compute numerically the solution $y^{\flat}(t)$ to the initial 
value problem of the approximated ODE,
\begin{eqnarray} \label{eq:recovered_ode_system}
   y' = \sum_{\ell=0}^{L} \, c_{\ell} \, \varphi_{\ell}( y ) \, , \quad 
   y(t_0) = \left( \, y_{1,1},\ldots,y_{n,1} \, \right)^{T} , \\
   c_{\ell} = \left( \, c_{\ell}^{(1)},\ldots,c_{\ell}^{(n)} \, \right)^{T}
\end{eqnarray}
provided all sample time points $t_{k,j}$ start exact synchronously. For
this numerical task any suitable integration scheme can be applied, i.e., 
in most cases, a decent one-step Runge-Kutta scheme or, alternatively,
a linearly implicit scheme (e.g. LIMEX \citep{ehrig1999advanced,limex}) will do.

\medskip

If the discrepancy between the approximated solution $y^{\flat}(\cdot)$ 
and the given data $y_{k}(t_{k,j})$ is satisfactory we could say the 
unknown model has successfully been verified. 

\medskip

Or else, if the discrepancy is too high, we could invoke yet another
parameter identification run: This time with the recovered ODE system
(\ref{eq:recovered_ode_system}), and the unknown coefficients $c_{\ell}$ 
as parameters with an obvious starting guess.

\medskip

If all else fails, another ansatz function family has to be chosen, or 
the filtering parameters have to be adjusted differently, or the problem
at hand needs just more (measurement) data.

\section{Numerical Results and Discussion}

In this section, we will study a predator-prey model and a damped pendulum model.
For the ansatz function family we choose in both cases multi-variate polynomials,
\begin{equation}
   \varphi_{\ell}(x) := \varphi_{(\ell_1,\ldots,\ell_n)}(x_{1},\ldots,x_{n}) := 
      x_{1}^{\ell_{1}} \cdot x_{2}^{\ell_{2}} \cdot \ldots \cdot x_{n}^{\ell_n}  
\end{equation}
where $\ell = (\ell_1,\ldots,\ell_n)$ denotes a multi-index with the usual 
meaning, as indicated. Consequently, for a polynomial of maximal degree $d$, 
i.e.~$|\ell| = \ell_1 + \ldots + \ell_n \leq d$, there are
\[
    \left( n+d \atop d \right) = \left( n+d \atop n \right)
\]
different monomial terms $\varphi_{\ell}(\cdot)$ with the determining 
coefficients $c_{\ell} \in \mathbb{R}$, specifying uniquely the Gram 
matrix $A$ in (\ref{eq:lsq_linear_system}).

\medskip

The runtime protocols for each component of the example computations, that will 
be given below, have the following structure.  In the first and second column, 
two different representations of the multi-index are listed (of which the one in the 
second column is the more familiar notation).  In the third column the percentage 
taken with respect to the absolute maximum of the coefficients $c_{\ell}$ found 
is given, and finally, in the last column, the actual value of each coefficient 
is printed.  Additionally, the resulting equation for the investigated component 
is recorded in which only monomial terms are included according to the indicated 
percentage threshold.  Last but not least, the weighted $\ell_2$-norm of the residual 
of the least-squares equation (\ref{eq:lsq_linear_system}) is displayed.

\subsection{Hare and Lynx}
\label{sec:hare_lynx}

Commonly, the dynamics of a predator-prey system is modelled by Lotka-Volterra-type 
equations \citep{Dd2015_GUIDE},
\begin{eqnarray}
   \label{eq:hare_lynx_ode1}
   y_{\mathrm{prey}}' & =   & y_{\mathrm{prey}} \, (\alpha - \beta \, y_{\mathrm{pred}}) \, , \\
   \label{eq:hare_lynx_ode2}
   y_{\mathrm{pred}}' & = - & y_{\mathrm{pred}} \, (\gamma - \delta \, y_{\mathrm{prey}})
\end{eqnarray}
where $y_{\mathrm{prey}}$ and $y_{\mathrm{pred}}$ denote the population 
number of the prey and the predator, respectively.  Here, the coupling 
constants $\alpha,\beta,\gamma,\delta \in \mathbb{R}_{\geq 0}$ describe the influence 
of the species on each other.  We want to recover such a system from some 
historical data of numbers of lynxes (predator) and hares (prey), collected 
by the Hudson Bay Company in Canada during the years $1900$--$1920$.  This 
data set can be seen, more or less, as representative for the total population 
number of both species \citep{Dd2015_GUIDE}.

\medskip

For basic biological reasons a constant term is missing in both components of the 
anticipated predator-prey model (\ref{eq:hare_lynx_ode1})--(\ref{eq:hare_lynx_ode2}).
Hence, we exclude the constant monomial term in both reconstruction runs as well.  
The maximal degree of the multi-variate polynomials is set to 3.

\begin{figure}[ht]

   \begin{center}
           {\bf (a)} \includegraphics[scale=.33]{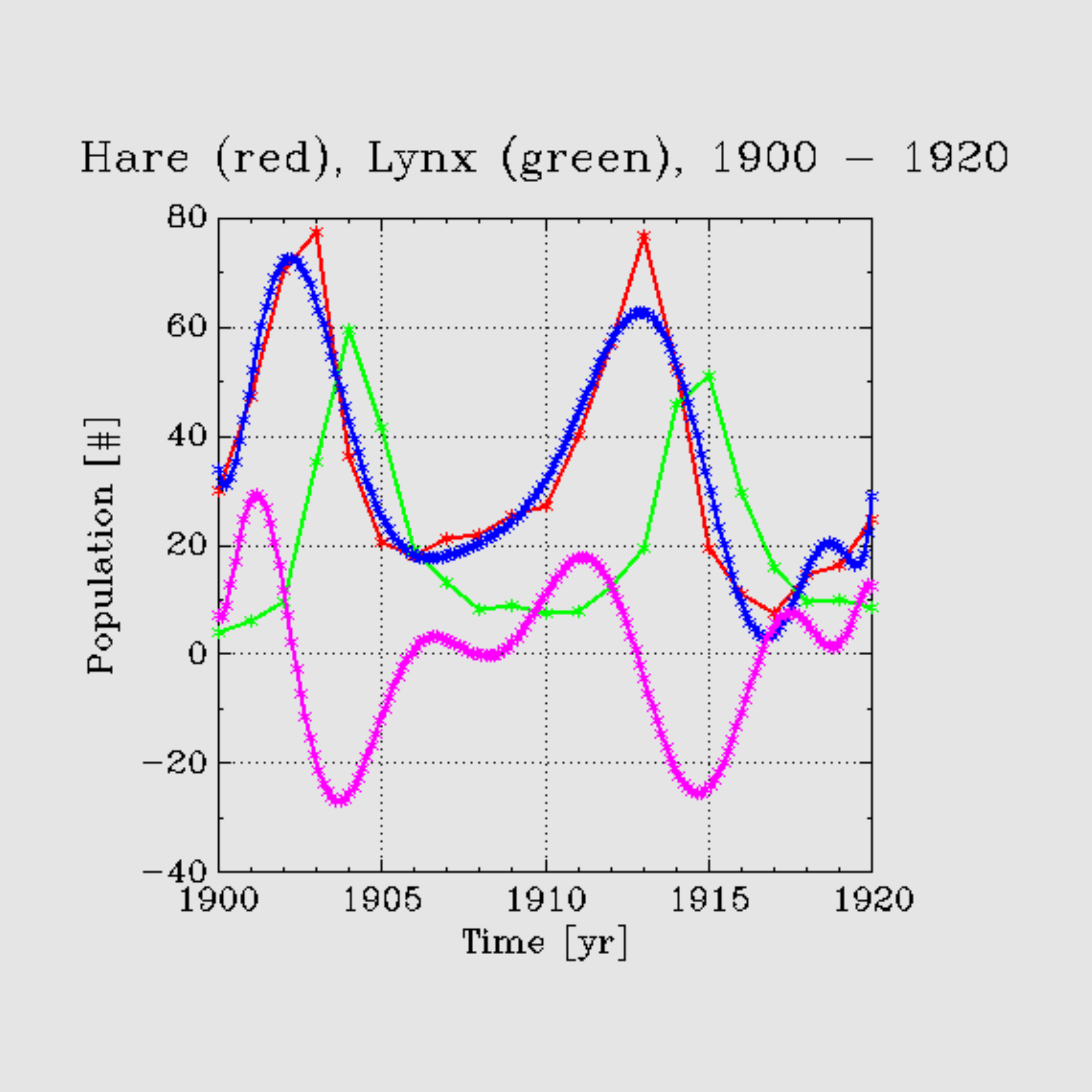} 
       \ \ {\bf (b)} \includegraphics[scale=.33]{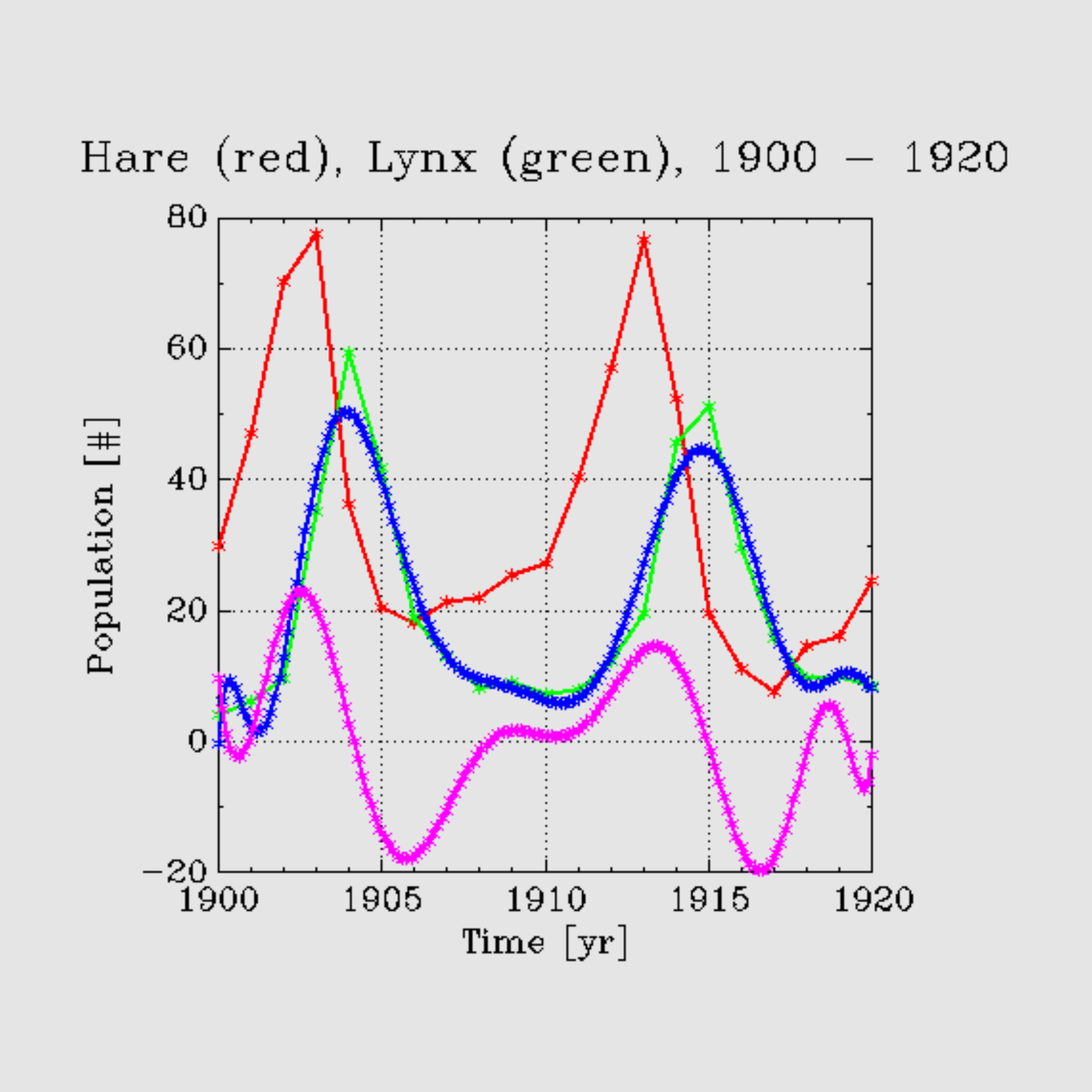}
      \caption{$T_n(x)$ approximations (blue) and its derivative (magenta) to 
               {\bf (a)} hare data, and {\bf (b)} lynx data 
               (n=80, evaluation of the approximations with n=11 each).} 
   \end{center}

\end{figure}

\vspace{-5ex}

\paragraph{Runtime protocol for $y_{0} = y_{\mathrm{prey}}$ component.}
\begin{Verbatim}[frame=single]
#total = 16   (max. deg. 3)
 [ 0 2 ] --> [0, 1]   100.00 %   c_k =  8.999962e-01
 [ 0 3 ] --> [0, 2]     3.79 %   c_k = -3.408105e-02
 [ 0 4 ] --> [0, 3]     0.00 %   c_k =  1.703492e-04
 [ 1 2 ] --> [1, 0]    19.66 %   c_k = -1.769268e-01
 [ 1 3 ] --> [1, 1]     3.44 %   c_k = -3.092972e-02
 [ 1 4 ] --> [1, 2]     0.00 %   c_k =  3.976574e-04
 [ 2 3 ] --> [2, 0]     2.65 %   c_k =  2.384418e-02
 [ 2 4 ] --> [2, 1]     0.00 %   c_k = -4.353977e-05
 [ 3 4 ] --> [3, 0]     0.00 %   c_k = -2.065367e-04
m =  9 monomial(s)      0.10 %   threshold
f(y0,y1) = + ( 9.0e-01) y1^1 + (-3.4e-02) y1^2 + (-1.8e-01) y0^1 
           + (-3.1e-02) y0^1 y1^1 + ( 2.4e-02) y0^2 
LSQ: ||residual/sqrt(n)||_2 = 3.095545072754599
\end{Verbatim}

\medskip

The $y_{\mathrm{prey}}$ is found to be described by the equation
\begin{equation}
   y_{0}' = + \, y_{0} \, ((a \, y_{0} - b) - c \, y_{1}) + d \, y_{1} - e \, y_{1}^2
\end{equation}
with $a = 0.024$, $b = 0.18$, $c = 0.031$, $d = 0.9$, and $e = 0.034$ .

\newpage 


\begin{figure}[ht]

   \begin{center}
           {\bf (a)} \includegraphics[scale=.33]{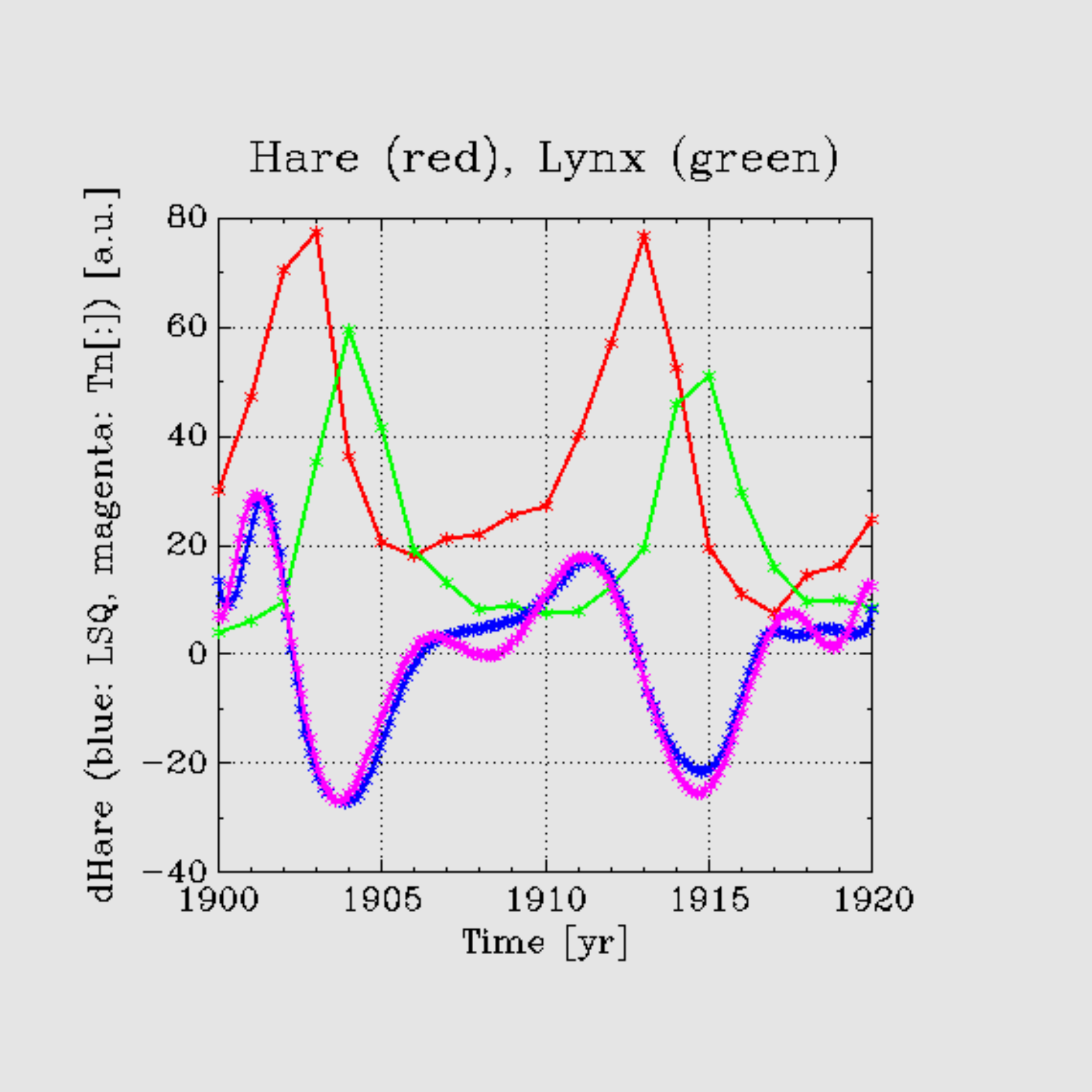} 
       \ \ {\bf (b)} \includegraphics[scale=.33]{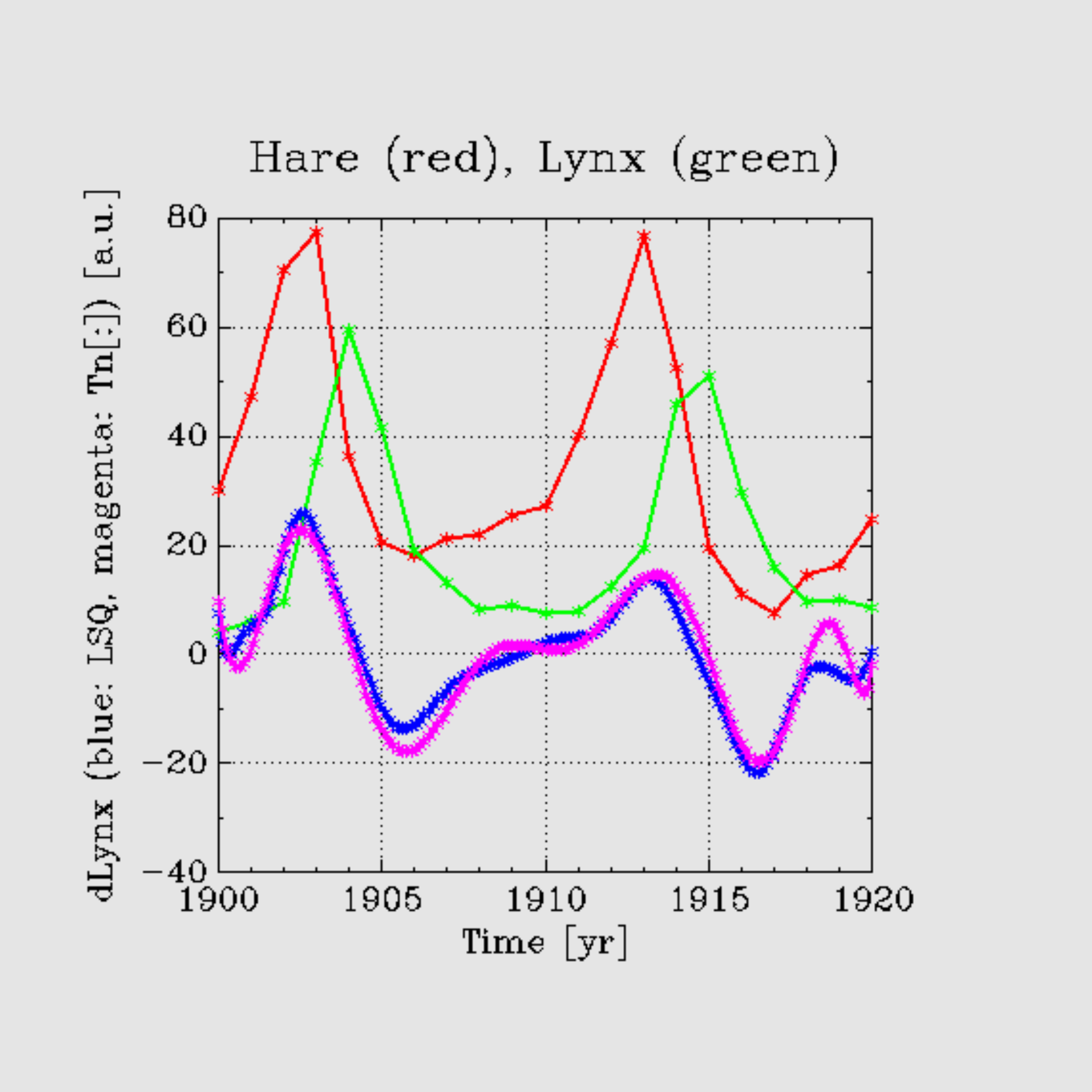}
      \caption{\label{fig:lsq_solution}
               Least-squares solution (blue) with $(c^{(k)})^*, k=1,2$, and 
               degree $=3$ to the unknown $F:y^* \mapsto (y^*)'$ for 
               {\bf (a)} hare data, and {\bf (b)} lynx data 
               (magenta each, and evaluation with 150 equidistant time point each).} 
   \end{center}
   
\end{figure}

\paragraph{Runtime protocol for $y_{1} = y_{\mathrm{pred}}$ component.}
\begin{Verbatim}[frame=single]
#total = 16   (max. deg. 3)
 [ 0 2 ] --> [0, 1]   100.00 %   c_k = -1.590145e+00
 [ 0 3 ] --> [0, 2]     1.35 %   c_k =  2.146557e-02
 [ 0 4 ] --> [0, 3]     0.00 %   c_k =  1.283903e-04
 [ 1 2 ] --> [1, 0]    45.70 %   c_k =  7.266958e-01
 [ 1 3 ] --> [1, 1]     0.87 %   c_k =  1.379691e-02
 [ 1 4 ] --> [1, 2]     0.00 %   c_k = -4.190196e-04
 [ 2 3 ] --> [2, 0]     1.35 %   c_k = -2.140643e-02
 [ 2 4 ] --> [2, 1]     0.00 %   c_k =  2.619544e-04
 [ 3 4 ] --> [3, 0]     0.00 %   c_k =  1.834174e-04
m =  9 monomial(s)      0.10 %   threshold
f(y0,y1) = + (-1.6e+00) y1^1 + ( 2.1e-02) y1^2 + ( 7.3e-01) y0^1 
           + ( 1.4e-02) y0^1 y1^1 + (-2.1e-02) y0^2 
LSQ: ||residual/sqrt(n)||_2 = 3.1046555639780498
\end{Verbatim}

\medskip

Here, the $y_{\mathrm{pred}}$ component is found to be given by the equation
\begin{equation}
   y_{1}' = - y_{1} \, (a  - b \, y_{0}) + c \, y_{1}^2 + d \, y_{0} - e \, y_{0}^2
\end{equation}
with $a = 1.6$, $b = 0.014$, $c = 0.021$, $d = 0.73$, and $e = 0.021$ .

\newpage

It would be interesting to see if an additional, iterative Gauss-Newton fit
would improve the identified parameters further and, possibly, establish some
link to the standard predator-prey model (\ref{eq:hare_lynx_ode1})--(\ref{eq:hare_lynx_ode2}).
Unfortunately, for the case with multi-variate polynomials of maximal degree 3,
as presented here, such an additional verification has not been successful at all.
Instead, the Gauss-Newton iteration stops after a few steps with the ODE system being
not integrable any more, i.e.~the Newton path of the coefficient sets during the 
iteration inexplicably leads to an ODE system that can not be solved numerically.
A systematic and comprehensive investigation of this unexpected behaviour would 
certainly go beyond the scope of this article and hence, must be left open for now.

\medskip

However, a repeated computation, completely analogous to the one presented here, 
only with the maximal polynomial degree restricted to 2, nicely demonstrates a 
relationship to the standard predator-prey model, as we have previously gathered,
see Table \ref{tab:coeffs_comparison}.  In particular, there is an intriguing
match between the coefficients of the second component, i.e.~the ODE equation
of the lynx.  The incompatibility factor in the last row of Table 
\ref{tab:coeffs_comparison} is determined by the ratio of the norms of successive 
updates during the Gauss-Newton iteration and thus, can be seen as an estimate
of the convergence speed, especially at the end of the iteration.  This factor,
if strictly below $1.0$, can be thought to measure how compatible a model to 
given data is (see \citep[sec.~4.3.2]{Dd_newton} for the mathematical background).  
The values in the last column of Table \ref{tab:coeffs_comparison} can be found 
in \citep{Dd2015_GUIDE}.  All the computational details are given in the Appendix.

\begin{table}[ht]
 \caption{\label{tab:coeffs_comparison} 
          Identified coefficients of the standard model \citep{Dd2015_GUIDE} and our data-based approach.}
 \begin{indented}
 \item[]\begin{tabular}{@{}lcrrl}
   \br
   Component                         & Multi-index &     Data-based  &      Standard  &             \\
   \mr
                                     &   $[0,1]$   &   1.848730e-01  &                &             \\ 
                                     &   $[0,2]$   &   3.675366e-03  &                &             \\
   $y_0'=f_{\mathrm{prey}}(y_0,y_1)$ &   $[1,0]$   &   2.531473e-01  &  5.475337e-01  & ($=\alpha$) \\
                                     &   $[1,1]$   &  -3.789005e-02  & -2.811932e-02  & ($=-\beta$) \\
                                     &   $[2,0]$   &   7.856857e-03  &                &             \\
                                     &             &                 &                &             \\
                                     &   $[0,1]$   &  -9.587181e-01  & -8.431750e-01  & ($=-\gamma$)\\ 
                                     &   $[0,2]$   &   6.845210e-04  &                &             \\
   $y_1'=f_{\mathrm{pred}}(y_0,y_1)$ &   $[1,0]$   &   5.296776e-02  &                &             \\
                                     &   $[1,1]$   &   2.408610e-02  &  2.655759e-02  & ($=\delta$) \\
                                     &   $[2,0]$   &  -1.113486e-04  &                &             \\
                                     &             &                 &                &             \\
   Scaled residual ({\tt Normf})     &             &     3.0712e-00  &    4.2362e-00  &             \\
   $\kappa$ (incomp. factor, see \citep{Dd_newton})         
                                     &             &     1.0274e-01  &    2.9798e-02  &             \\
   \br
  \end{tabular}
  \end{indented}
\end{table}

\medskip

\begin{figure}[ht]
   \begin{center}
       \ \ \ \ \includegraphics[scale=.38]{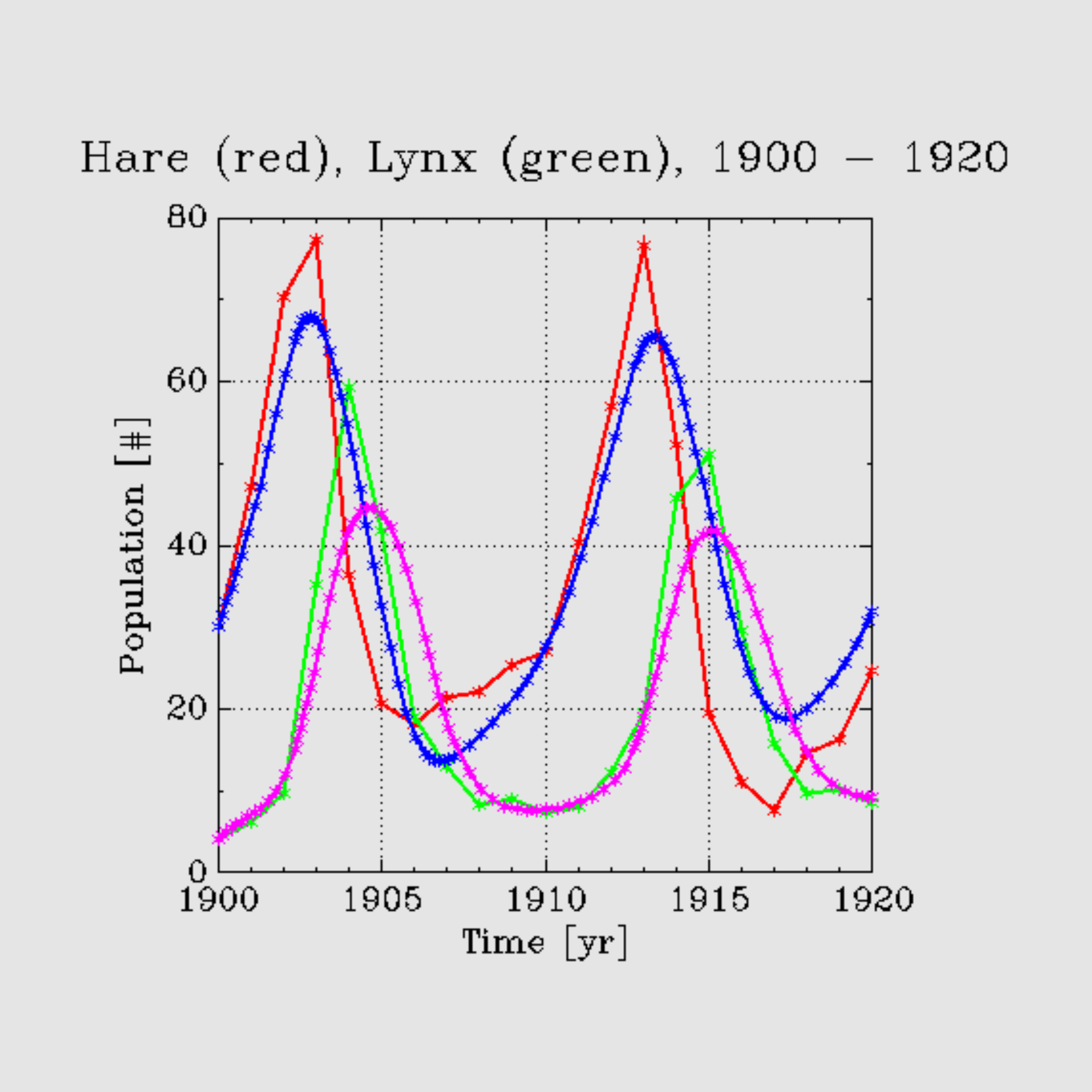} 
      \caption{\label{fig:hare_lynx_verification} 
               Comparison between given data (red, green) and solution (blue, magenta) to
               the approximated ODE with multi-variate polynomials of maximal degree $3$,
               as shown in Fig.~\ref{fig:lsq_solution}
               (by using the adaptive time stepping integrator LIMEX).} 
   \end{center}

\end{figure}

To sum up the hare and lynx case, the reconstructed components $y_0$ and $y_1$ seem 
to follow only partly the standard predator-prey model.  Instead, our findings here
show that the data might be more involved than previously assumed.  Our conclusion
is also supported by the verification result, as can be seen in Figure 
\ref{fig:hare_lynx_verification}.


\subsection{Damped Pendulum}

Here, we consider the initial value problem
\begin{eqnarray}
   \theta_1' & = \theta_2 \, ,  \label{eq:pendulum_ode1} \\
   \theta_2' & = -u \, \theta_2 - \frac{g}{l} \, \sin \theta_1 \label{eq:pendulum_ode2}
\end{eqnarray}
with the initial condition $(\theta_1, \theta_2)(0) = (1,0)$, and the constants 
$u = 0.25$, $l = 2.0$, and $g = 9.81$ (the gravitational acceleration on earth) 
are used \citep{perona2000trajectory}.  Very sparse data has been generated by solving 
numerically (\ref{eq:pendulum_ode1})--(\ref{eq:pendulum_ode2}) in the time interval 
$[0, 10]$, and subsequently taking $49$ equidistant samples of the solution 
trajectories in this interval.

\medskip

In order to be able to catch a glimpse of the non-linearity on the right-hand side 
in (\ref{eq:pendulum_ode2}), for the reconstruction trial the maximal degree of the 
multi-variate polynomials is set to $4$.

\begin{figure}[ht]

   \begin{center}
           {\bf (a)} \includegraphics[scale=.33]{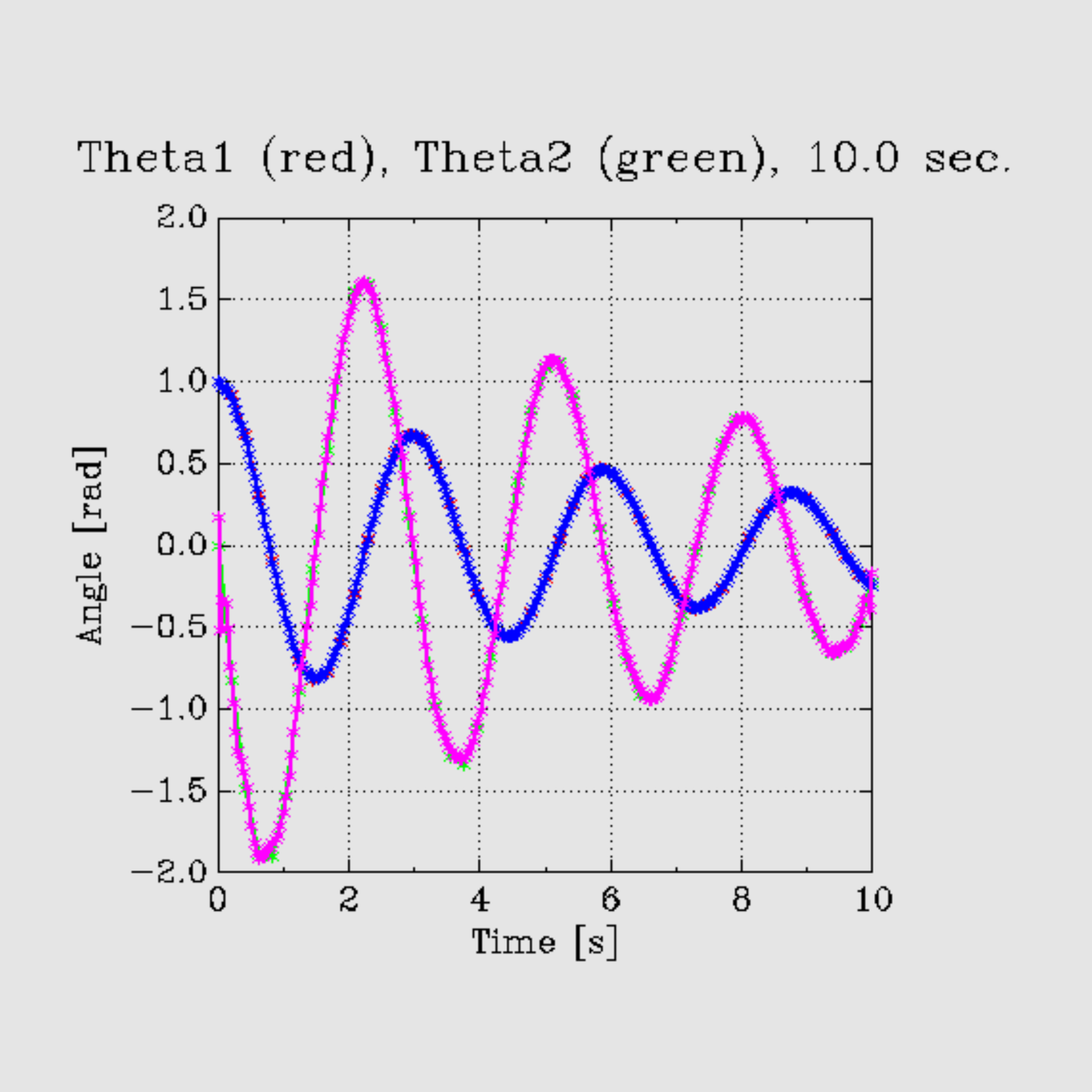} 
       \ \ {\bf (b)} \includegraphics[scale=.33]{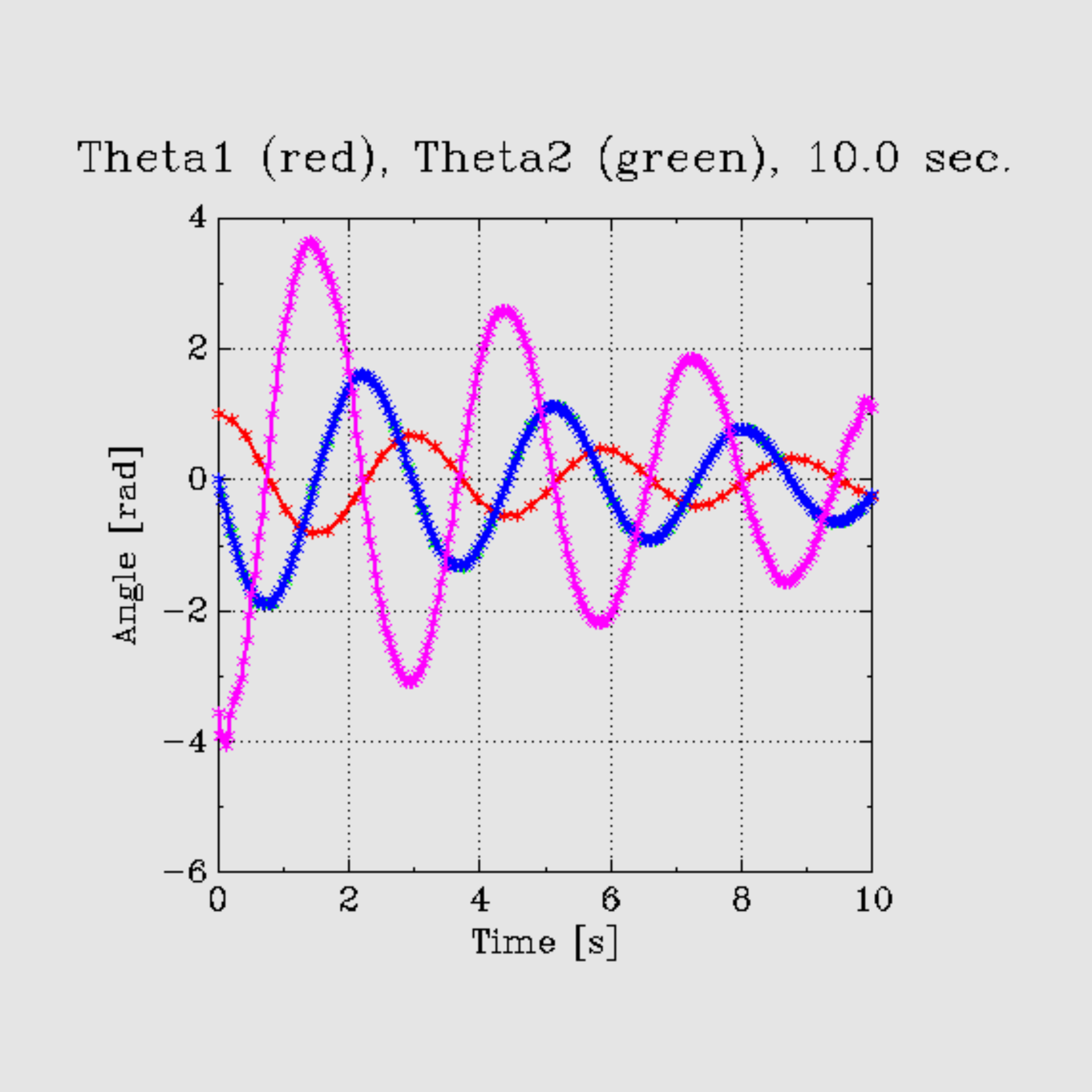}
      \caption{$T_n(x)$ approximations (blue) and its derivative (magenta) to 
               {\bf (a)} $\theta_1$ data, and {\bf (b)} $\theta_2$ data 
               (n=80, evaluation of the approximations with n=62 each).} 
   \end{center}

\end{figure}

\medskip

As seen in the following runtime protocol for the $\theta_{1}$ component, 
the equation
\begin{equation}
  y_{0}' = a \, y_{1}
\end{equation}
with $a = 1.0$ is found.

\newpage

\paragraph{Runtime protocol for $y_{0} = \theta_{1}$ component.}
\begin{Verbatim}[frame=single]
#total = 25   (max. deg. 4)
 [ 0 1 ] --> [0, 0]     0.00 %   c_k =  2.001792e-03
 [ 0 2 ] --> [0, 1]   100.00 %   c_k =  9.982190e-01
 [ 0 3 ] --> [0, 2]     0.00 %   c_k = -6.094577e-03
 [ 0 4 ] --> [0, 3]     0.00 %   c_k = -2.859434e-04
 [ 0 5 ] --> [0, 4]     0.00 %   c_k =  1.027003e-03
 [ 1 2 ] --> [1, 0]     0.00 %   c_k =  1.384047e-03
 [ 1 3 ] --> [1, 1]     0.00 %   c_k =  2.514082e-02
 [ 1 4 ] --> [1, 2]     0.00 %   c_k =  4.947748e-03
 [ 1 5 ] --> [1, 3]     0.00 %   c_k = -5.810996e-03
 [ 2 3 ] --> [2, 0]     0.00 %   c_k = -2.636188e-02
 [ 2 4 ] --> [2, 1]     0.00 %   c_k =  3.051582e-02
 [ 2 5 ] --> [2, 2]     0.00 %   c_k =  2.399808e-02
 [ 3 4 ] --> [3, 0]     0.00 %   c_k = -2.561132e-03
 [ 3 5 ] --> [3, 1]     0.00 %   c_k = -1.691789e-02
 [ 4 5 ] --> [4, 0]     0.00 %   c_k =  3.776657e-02
m = 15 monomial(s)      5.00 %   threshold
f(y0,y1) = + ( 1.0e+00) y1^1 
LSQ: ||residual/sqrt(n)||_2 = 0.035592552825152
\end{Verbatim}


\newpage

\begin{figure}[ht] 
   \begin{center}
           {\bf (a)} \includegraphics[scale=.33]{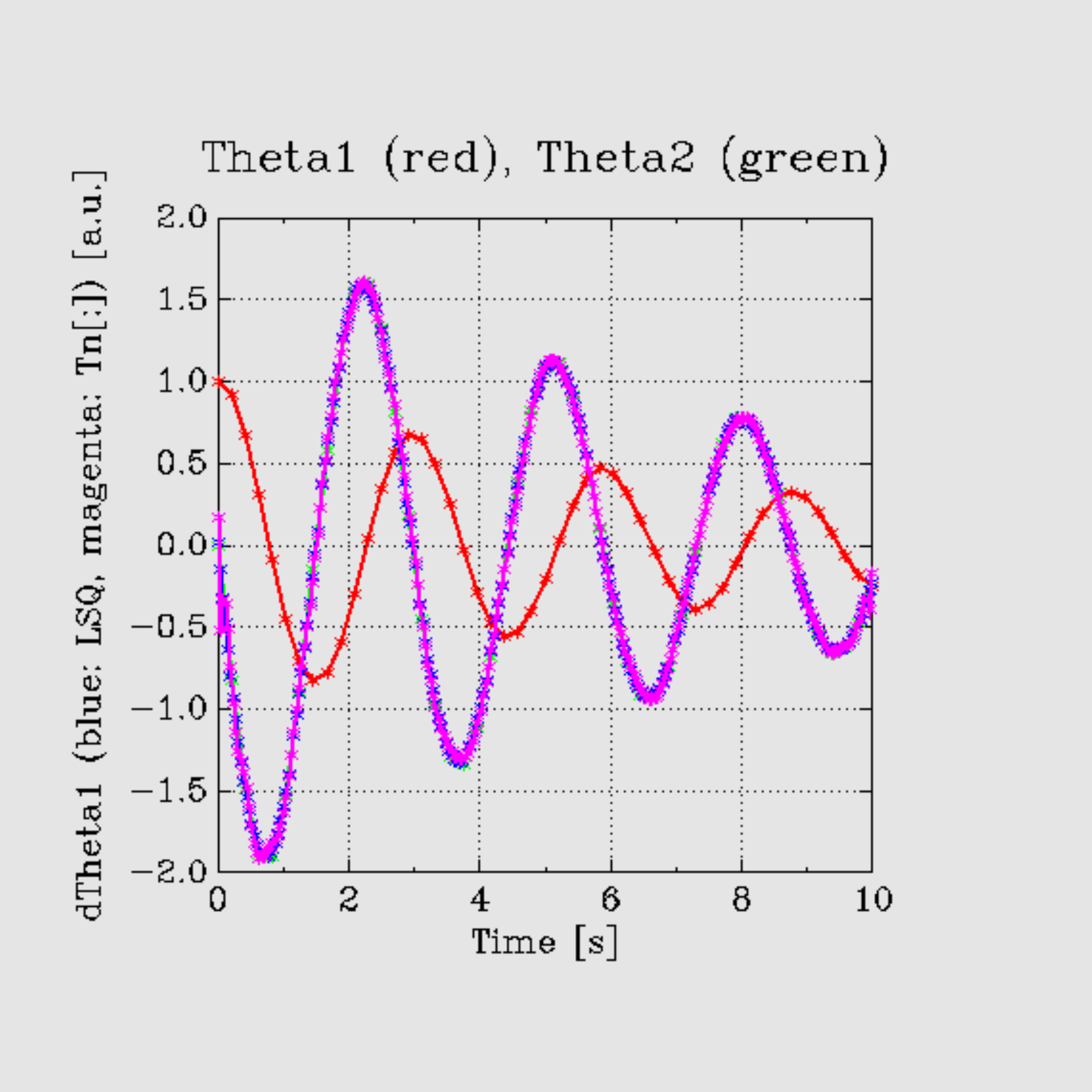} 
       \ \ {\bf (b)} \includegraphics[scale=.33]{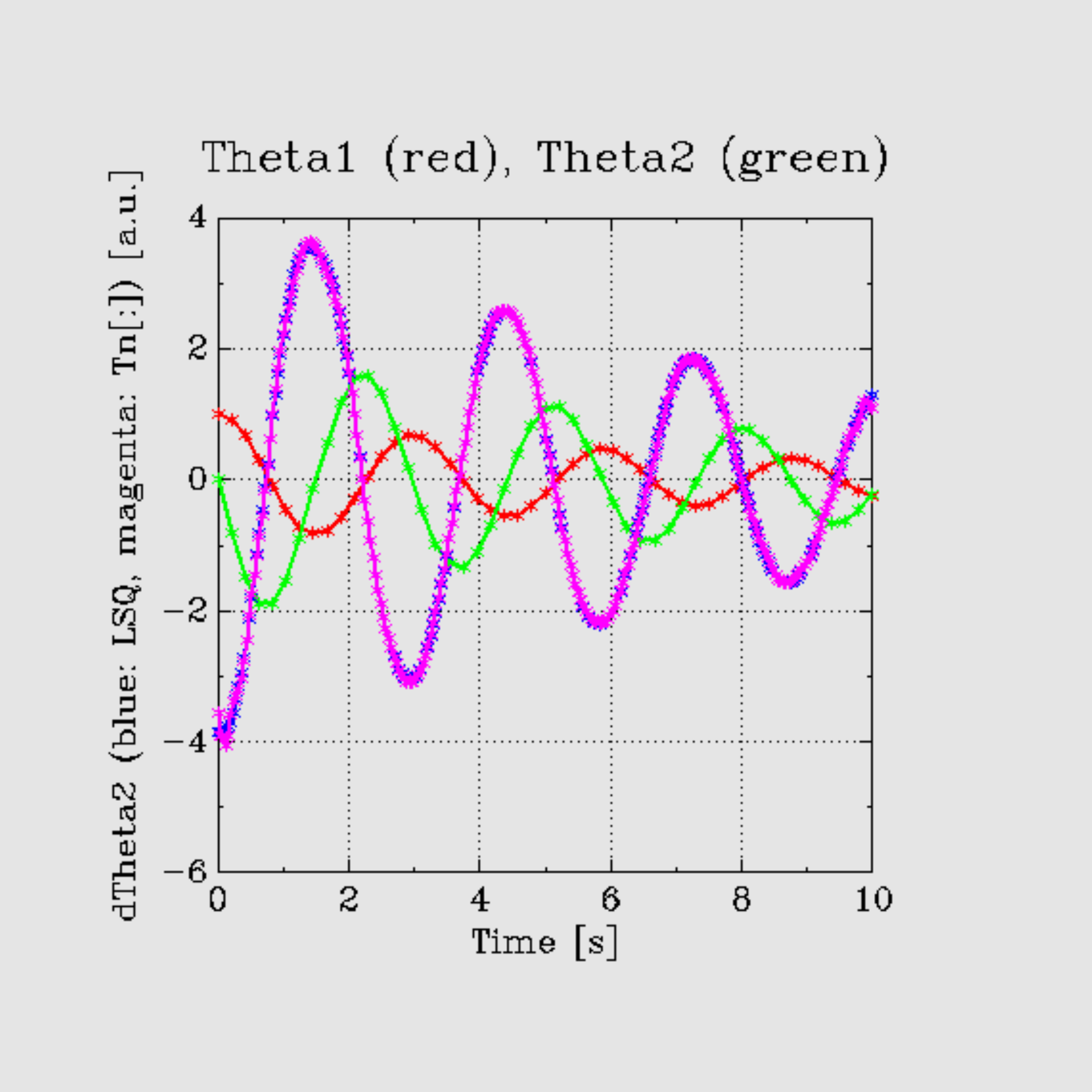}
      \caption{\label{fig:lsq_solution_dp}
               Least-squares solution (blue) with $(c^{(k)})^*, k=1,2$, and 
               degree $=4$ to the unknown $F:y^* \mapsto (y^*)'$ for 
               {\bf (a)} $\theta_1$ data, and {\bf (b)} $\theta_2$ data 
               (magenta each, and evaluation with 250 equidistant time point each).} 
   \end{center}
   
\end{figure}


\paragraph{Runtime protocol for $y_{1} = \theta_{2}$ component.}
\begin{Verbatim}[frame=single]
#total = 25   (max. deg. 4)
 [ 0 1 ] --> [0, 0]     0.00 %   c_k =  5.555385e-03
 [ 0 2 ] --> [0, 1]     5.13 %   c_k = -2.535032e-01
 [ 0 3 ] --> [0, 2]     0.00 %   c_k = -1.796957e-02
 [ 0 4 ] --> [0, 3]     0.00 %   c_k = -6.444822e-03
 [ 0 5 ] --> [0, 4]     0.00 %   c_k = -3.850649e-04
 [ 1 2 ] --> [1, 0]   100.00 %   c_k = -4.942789e+00
 [ 1 3 ] --> [1, 1]     0.00 %   c_k = -1.564249e-01
 [ 1 4 ] --> [1, 2]     0.00 %   c_k =  2.928157e-02
 [ 1 5 ] --> [1, 3]     0.00 %   c_k =  1.852611e-02
 [ 2 3 ] --> [2, 0]     0.00 %   c_k = -9.883681e-02
 [ 2 4 ] --> [2, 1]     0.00 %   c_k =  1.441293e-01
 [ 2 5 ] --> [2, 2]     0.00 %   c_k =  1.787311e-01
 [ 3 4 ] --> [3, 0]    19.50 %   c_k =  9.639487e-01
 [ 3 5 ] --> [3, 1]     9.90 %   c_k =  4.893068e-01
 [ 4 5 ] --> [4, 0]     0.00 %   c_k =  2.128485e-01
m = 15 monomial(s)      5.00 %   threshold
f(y0,y1) = + (-2.5e-01) y1^1 + (-4.9e+00) y0^1 + ( 9.6e-01) y0^3 
           + ( 4.9e-01) y0^3 y1^1 
LSQ: ||residual/sqrt(n)||_2 = 0.04531770760243261
\end{Verbatim}

\newpage

In this case, for the $\theta_2$ component, we obtain
\begin{equation}
   y_{1}' = - a \, y_{1} - b \, y_{0} + c \, y_{0}^3 + d \, y_{0}^3 \, y_{1} 
\end{equation}
with $a = 0.25$, $b = 4.9$, $c = 0.96$, and $d = 0.49$ .

\medskip

\begin{figure}[ht]
   \begin{center}
       \ \ \ \ \includegraphics[scale=.38]{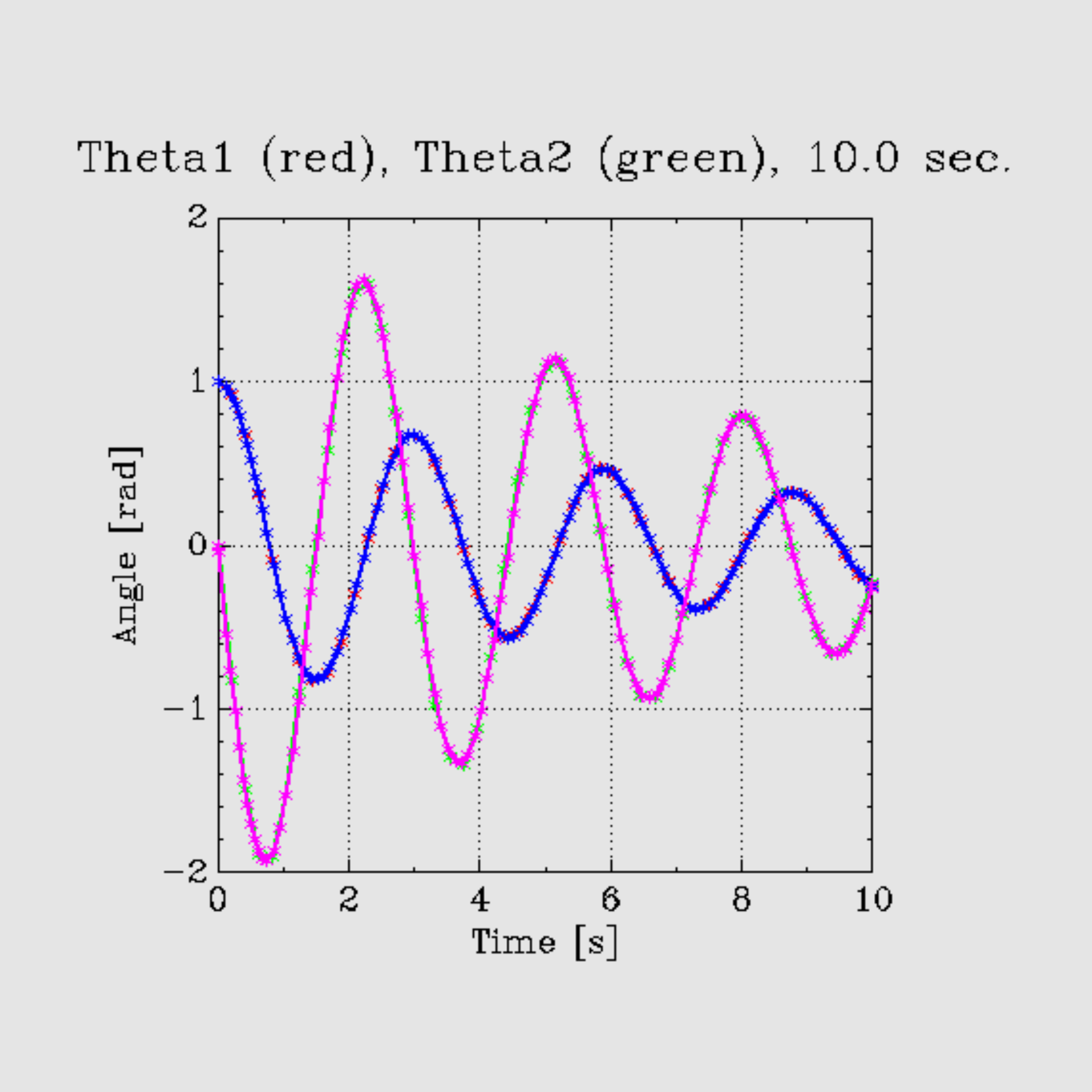} 
      \caption{\label{fig:damped_pendulum_verification}
               Comparison between given data (red, green) and solution (blue, magenta) to
               the approximated ODE with multi-variate polynomials of maximal degree $4$,
               as shown in Fig.~\ref{fig:lsq_solution_dp}
               (by using the adaptive time stepping integrator LIMEX).} 
   \end{center}

\end{figure}

Summarising these findings, the recovered equations in this non-linear case 
agree convincingly well with the pendulum model that has been used to generated 
the given sparse data set.  Especially for the second component $y_1=\theta_2$, 
the reconstructed coefficients $b = 4.943 \sim g/l = 4.905$ and 
$c = 0.9639 \sim -g/6l = 0.8175$ match the Taylor coefficients of the sine in 
the right-hand side of the original damped pendulum model.  The discrepancy 
in $c$ is compensated by the additional coefficient $d$, indicating that the 
approximated model is still not complete.  Nevertheless, the verification result, 
as can be seen clearly in Figure \ref{fig:damped_pendulum_verification}, is already
almost perfect.

\section{Conclusion}

In the present paper, a highly flexible method for the reconstruction of an 
unknown dynamical system, to be described in terms of an explicit ODE 
system, is presented.  The reconstruction method is based on sampled, 
possibly sparse, trajectory solution points as given data.

\medskip

In principle, any kind of relationship between the available information
is tried to be recovered automatically, if different components of the 
unknown system are given.  Yet, because of the extreme simple structure 
of the presented method, and the fact that the ODE itself is not being
solved during the reconstruction, it is possible to readily adapt the 
approach to any anticipated behaviour of the unknown dynamical system.
Moreover, the approach is fast and, most importantly, highly reliable.

\medskip

Additionally, the results of our fast reconstruction method could be 
used as a very good starting guess, for example, of more sophisticated 
iterative identification procedures such as well-established Gauss-Newton 
codes \citep{Dd_newton}.


%
%
%
%

\ack
The author would like to thank Sebastian G\"otschel for his invaluable comments and 
clarifications.  Additionally, the author also wishes to thank the CSB group at ZIB 
for many fruitful and lively discussions.

\newpage

\appendix
\section*{Appendix}
\setcounter{section}{1}

Here, we repeat the computations of section \ref{sec:hare_lynx} for the predator-prey case, 
only the maximal degree of the multi-variate polynomials is set to 2.  This results in the 
following two runtime protocols.


\paragraph{Runtime protocol for $y_{0} = y_{\mathrm{prey}}$ component.}
\begin{Verbatim}[frame=single]
#total = 9   (max. deg. 2)
 [ 0 2 ] --> [0, 1]     7.72 %   c_k = -3.591668e-02
 [ 0 3 ] --> [0, 2]     0.93 %   c_k = -4.319190e-03
 [ 1 2 ] --> [1, 0]   100.00 %   c_k =  4.653563e-01
 [ 1 3 ] --> [1, 1]     3.58 %   c_k = -1.665442e-02
 [ 2 3 ] --> [2, 0]     0.00 %   c_k = -9.424325e-05
m =  5 monomial(s)      0.10 %   threshold
f(y0,y1) = + (-3.6e-02) y1^1 + (-4.3e-03) y1^2 + ( 4.7e-01) y0^1 
           + (-1.7e-02) y0^1 y1^1 
LSQ: ||residual/sqrt(n)||_2 = 3.7658055884421944
\end{Verbatim}


\paragraph{Runtime protocol for $y_{1} = y_{\mathrm{pred}}$ component.}
\begin{Verbatim}[frame=single]
#total = 9   (max. deg. 2)
 [ 0 2 ] --> [0, 1]   100.00 %   c_k = -9.800874e-01
 [ 0 3 ] --> [0, 2]     0.68 %   c_k =  6.673957e-03
 [ 1 2 ] --> [1, 0]     7.26 %   c_k =  7.117392e-02
 [ 1 3 ] --> [1, 1]     1.52 %   c_k =  1.487619e-02
 [ 2 3 ] --> [2, 0]     0.16 %   c_k =  1.545017e-03
m =  5 monomial(s)      0.10 %   threshold
f(y0,y1) = + (-9.8e-01) y1^1 + ( 6.7e-03) y1^2 + ( 7.1e-02) y0^1 
           + ( 1.5e-02) y0^1 y1^1 + ( 1.5e-03) y0^2 
LSQ: ||residual/sqrt(n)||_2 = 3.666500740586521
\end{Verbatim}

\medskip

Taking these coefficients $c_{\ell}^{\mathrm{prey}}$ and $c_{\ell}^{\mathrm{pred}}$,
$|\ell| \leq 2$, as starting values for an non-linear least squares method, in order 
to identify parameters of the ODE model build by these two polynomials
\begin{eqnarray}
   \label{eq:apprx_hare_lynx_ode1}
   y_{0}' & = & \sum\limits_{|\ell| \leq 2} c_{\ell}^{\mathrm{prey}} \, y^{\ell} \, , \\
   \label{eq:apprx_hare_lynx_ode2}
   y_{1}' & = & \sum\limits_{|\ell| \leq 2} c_{\ell}^{\mathrm{pred}} \, y^{\ell}
\end{eqnarray}
with initial conditions $y_0(1900)=30$ and $y_1(1900)=4$ kept fixed, will result in a 
runtime protocol as follows, trying to match the given measurement data best in the 
least squares sense by successively varying these 10 coefficients.  Here, we choose 
an error-oriented Gauss-Newton scheme with adaptive trust region and rank strategy,
{\tt NLSCON}, as fitting routine.  Solutions to the IVP of 
(\ref{eq:apprx_hare_lynx_ode1})--(\ref{eq:apprx_hare_lynx_ode2}) 
are computed by the linearly-implicit extrapolation integrator {\tt LIMEX} 
\citep{ehrig1999advanced,limex} during the Gauss-Newton iteration, with 
$\mathrm{RTOL}=1.0\mathrm{E}-9$ and $\mathrm{ATOL}=1.0\mathrm{E}-9$.
All other relevant settings are included in the resulting runtime protocol
of the parameter identification task for the data-based model 
(\ref{eq:apprx_hare_lynx_ode1})--(\ref{eq:apprx_hare_lynx_ode2}).

\paragraph{Runtime protocol for an error-oriented Gauss-Newton scheme (NLSCON).}
\begin{Verbatim}[frame=single]

     N L S C O N  *****  V e r s i o n  2 . 3 . 3 ***

 Gauss-Newton-Method for the solution of nonlinear least squares
 problems

 Real    Workspace declared as 1742 is used up to 1578 ( 90.6 percent) 
 Integer Workspace declared as   62 is used up to   50 ( 80.6 percent)
 
 
 Number of parameters to be estimated (N) :   10
 Number of data to fitted, e.g. observations (MFIT) :   42
 Number of equality constraints (MCON) :    0

 Prescribed relative precision (PTOL) :   0.10D-02

 The Jacobian is supplied by numerical differentiation 
                                          (feedback strategy included)
 Automatic row scaling of the Jacobian is allowed

 Rank-1 updates are inhibited
 Problem is specified as being highly nonlinear
 Bounded damping strategy is off
 Maximum permitted number of iteration steps :     40

 
 Internal parameters:

 Starting value for damping factor FCSTART =  0.10D-01
 Minimum allowed damping factor FCMIN =  0.10D-01
 Rank-1 updates decision parameter SIGMA =  0.10D+04
 Initial Jacobian pseudo-rank IRANK =    10
 Maximum permitted subcondition COND =  0.45D+16

***********************************************************************

    It       Normf               Normx       Damp.Fct.   New      Rank
     0      0.8345461D+01       0.251D+00                0        10
     1      0.8277889D+01    *  0.248D+00      0.010
 1
     1      0.8277889D+01       0.241D+00                0        10
     2      0.5289067D+01    *  0.916D-01      0.375
 2
     2      0.5289067D+01       0.991D-01                0        10
     2      0.9455332D+01    *  0.133D+00      0.819
     3      0.4570706D+01    *  0.755D-01      0.232
 3
     3      0.4570706D+01       0.827D-01                0        10
     4      0.6342464D+01    *  0.732D-01      1.000
 4
     4      0.6342464D+01       0.229D-01                0        10
     5      0.3444394D+01    *  0.723D-02      1.000
 5
     5      0.3444394D+01       0.171D-01                0        10
     6      0.3083738D+01    *  0.315D-02      0.955
 6
     6      0.3083738D+01       0.652D-02                0        10
     7      0.3071453D+01    *  0.120D-03      1.000
 7
     7      0.3071453D+01       0.176D-02                0        10
     8      0.3071223D+01    *  0.252D-05      1.000
 8
     8      0.3071223D+01       0.180D-03                0        10
     9      0.3071220D+01    *  0.957D-07      1.000

     
 Solution of nonlinear least squares problem obtained 
 within   9 iteration steps

 Incompatibility factor kappa 0.103D+00
 
 Achieved relative accuracy  0.207D-04 

 \end{Verbatim}


\paragraph{Additional statistical analysis of the result of the identification task.}
\begin{Verbatim}[frame=single]

   Standard deviation of parameters
   --------------------------------
     No.  Estimate           sigma(X)
      1   0.185D+00   +/-   0.174D+00     =   94.25 %
      2   0.368D-02   +/-   0.463D-02     =  125.87 %
      3   0.253D+00   +/-   0.107D+00     =   42.40 %
      4  -0.379D-01   +/-   0.646D-02     =   17.04 %
      5   0.786D-02   +/-   0.271D-02     =   34.46 %
      6  -0.959D+00   +/-   0.198D+00     =   20.67 %
      7   0.685D-03   +/-   0.491D-02     =  716.79 %
      8   0.530D-01   +/-   0.725D-01     =  136.89 %
      9   0.241D-01   +/-   0.465D-02     =   19.29 %
     10  -0.111D-03   +/-   0.124D-02     = 1109.90 %
     
   Independent confidence intervals
   --------------------------------
   (on 95%-probability level using F-distribution F(alfa,1,m-n)= 4.15)

      1  ( -0.170D+00 ,  0.540D+00 )
      2  ( -0.575D-02 ,  0.131D-01 )
      3  (  0.345D-01 ,  0.472D+00 )
      4  ( -0.510D-01 , -0.247D-01 )
      5  (  0.234D-02 ,  0.134D-01 )
      6  ( -0.136D+01 , -0.555D+00 )
      7  ( -0.931D-02 ,  0.107D-01 )
      8  ( -0.947D-01 ,  0.201D+00 )
      9  (  0.146D-01 ,  0.336D-01 )
     10  ( -0.263D-02 ,  0.241D-02 )

     
   ******  Statistics * NLSCON   *******
   ***  Gauss-Newton iter.:       9  ***
   ***  Corrector steps   :       1  ***
   ***  Rejected rk-1 st. :       0  ***
   ***  Jacobian eval.    :      10  ***
   ***  Function eval.    :      12  ***
   ***  ...  for Jacobian :     100  ***
   *************************************
\end{Verbatim}

\medskip

\begin{figure}[ht]

   \begin{center}
           {\bf (a)} \includegraphics[scale=.55]{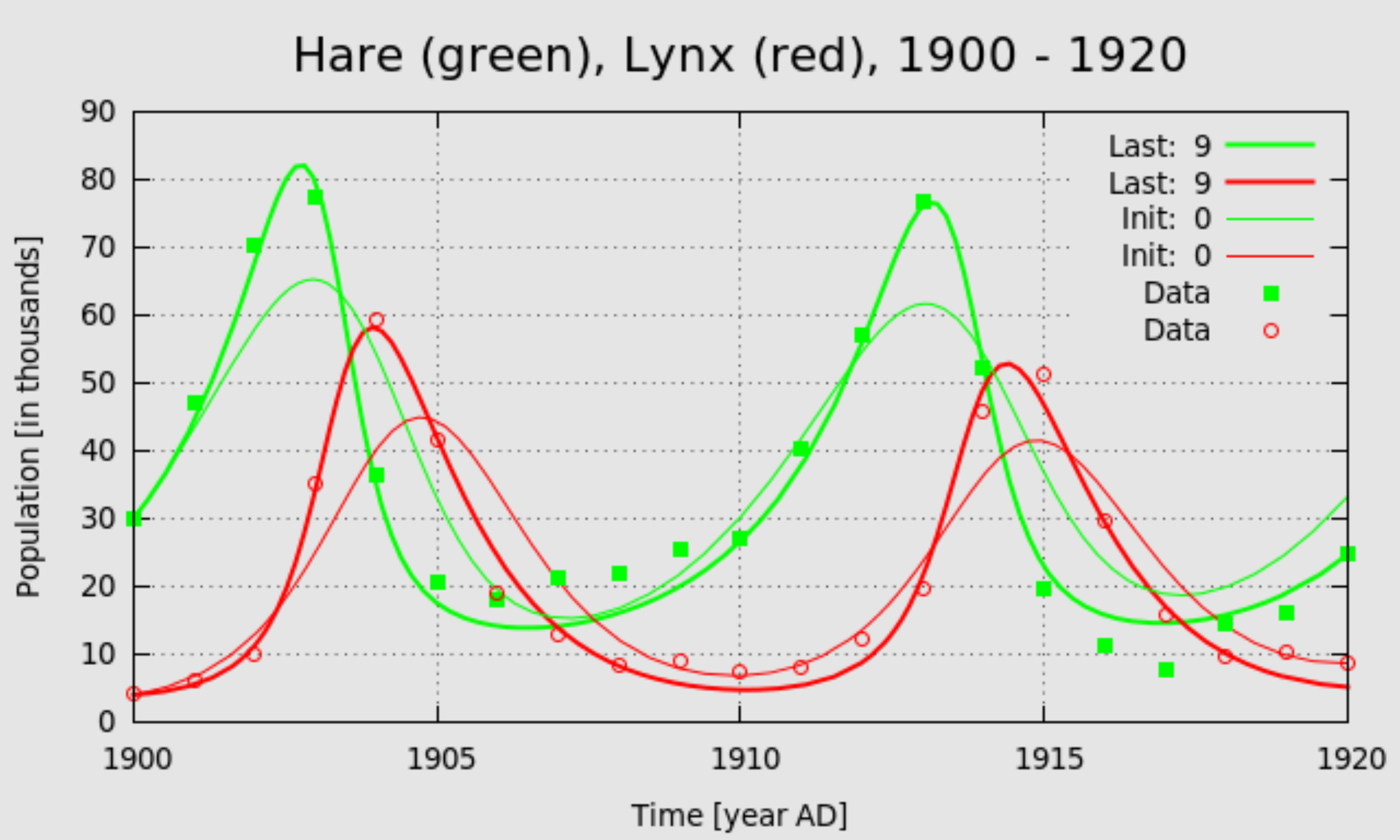} \\[3ex]
           {\bf (b)} \includegraphics[scale=.55]{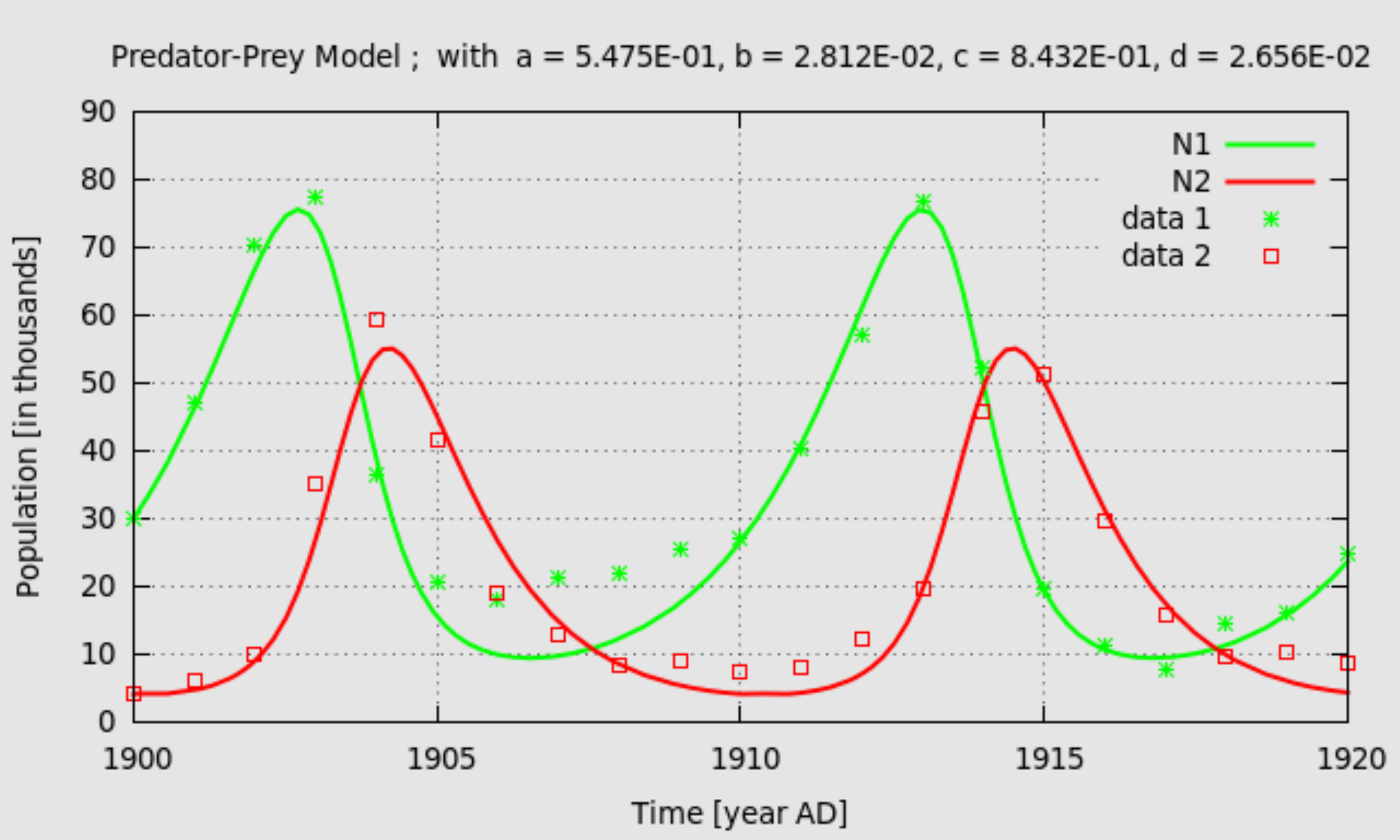}
      \caption{\label{fig:nlscon_solution}
               Fitted population curves {\bf (a)} of the data-based model (with maximal
               polynomial degree $=2$) and {\bf (b)} of the standard model \citep{Dd2015_GUIDE}.}
   \end{center}
   
\end{figure}

In order to have a visual inspection of these findings, the resulting population curves 
are plotted in Figure \ref{fig:nlscon_solution}.  In particular it seems that the curves
of our approximated model match more the given data in the time period 1905 -- 1915.

\medskip

Note that, although the Gauss-Newton iteration converges nicely, and with full rank of 
the Jacobians, the overall estimated standard deviation of the identified coefficients 
is relatively high.  This means that the linearisation at the minimum point found is 
nearly flat and hence, the corresponding confidence intervals are comparatively wide.

\bigskip

The implementation of the presented approach that, written in \texttt{Ruby}, 
has been applied to perform all example computations in this paper is available 
upon request.

%

\bibliographystyle{iopart-num} 
\bibliography{arXiv_GenericYetFlexible}

\providecommand{\newblock}{}
\begin{thebibliography}{10}
\expandafter\ifx\csname url\endcsname\relax
  \def\url#1{{\tt #1}}\fi
\expandafter\ifx\csname urlprefix\endcsname\relax\def\urlprefix{URL }\fi
\providecommand{\eprint}[2][]{\url{#2}}

\bibitem{Dd2015_GUIDE}
Deuflhard P and R\"oblitz S 2015 {\em {A Guide to Numerical Modelling in
  Systems Biology}\/} ({\em {Texts in Computational Science and Engineering}\/}
  no~12) ({Springer-Verlag})

\bibitem{perona2000trajectory}
Perona P, Porporato A and Ridolfi L 2000 {\em {Nonlinear Dynamics}\/} {\bf 23}
  13--33

\bibitem{walter1997identification}
Walter E and Pronzato L 1997 {\em {Identification of Parametric Models from
  Experimental Data}\/} (Springer)

\bibitem{Dd_NUM2}
Deuflhard P and Bornemann F 2002 {\em {Scientific Computing with Ordinary
  Differential Equations}\/} 1st ed ({\em {Texts in Applied Mathematics}\/}
  no~42) ({Springer-Verlag})

\bibitem{eisenhammer1991modeling}
Eisenhammer T, H\"ubler A, Packard N and Kelso J 1991 {\em {Biological
  Cybernetics}\/} {\bf 65} 107--112

\bibitem{raue2015data}
Raue A, Steiert B, Schelker M, Kreutz C, Maiwald T, Hass H, Vanlier J,
  T\"onsing C, Adlung L, Engesser R, Mader W, Heinemann T, Hasenauer J,
  Schilling M, H\"ofer T, Klipp E, Theis F, Klingm\"uller U, Sch\"oberl B and
  Timmer J 2015 {\em {Bioinformatics}\/} {\bf 31} 3558--3560

\bibitem{press2002_NR}
Press W~H, Teukolsky S~A, Vetterling W~T and Flannery B~P 2002 {\em {Numerical
  Recipes in C++ : The Art of Scientific Computing}\/} 2nd ed (Cambridge
  University Press)

\bibitem{trefethen2013_AT}
Trefethen L 2013 {\em {Approximation Theory and Approximation Practice}\/}
  (Philadelphia, PA: Society for Industrial and Applied Mathematics (SIAM))
  \urlprefix\url{http://www2.maths.ox.ac.uk/chebfun/ATAP/}

\bibitem{quarteroni2007_NM}
Quarteroni A, Sacco R and Saleri F 2007 {\em {Numerical Mathematics}\/} 2nd ed
  ({\em {Texts in Applied Mathematics}\/} no~37) ({Springer-Verlag})

\bibitem{de_boor1978_Guide}
De~Boor C 1978 {\em {A Practical Guide to Splines}\/} (New York:
  Springer-Verlag)

\bibitem{Dd_NUM1}
Deuflhard P and Hohmann A 2003 {\em {Numerical analysis in Modern Scientific
  Computing -- An Introduction}\/} 2nd ed ({\em {Texts in Applied
  Mathematics}\/} no~43) ({Springer-Verlag})

\bibitem{Dd_newton}
Deuflhard P 2004 {\em {Newton Methods for Nonlinear Problems -- Affine
  Invariance and Adaptive Algorithms}\/} ({\em {Springer Series in
  Computational Mathematics}\/} no~35) (Springer)

\bibitem{ehrig1999advanced}
Ehrig R, Nowak U, Oeverdieck L and Deuflhard P 1999 {\em {Lecture Notes in
  Computational Science and Engineering}\/} {\bf 8} 233--244

\bibitem{limex}
Deuflhard P and Nowak U 1987 Extrapolation integrators for quasilinear implicit
  {ODE}s {\em {Large Scale Scientific Computing}\/} ed Deuflhard P and Engquist
  B (Birkh\"auser) pp 37--50

\end{thebibliography}

\end{document}